# A Discrete Power Distribution

Subrata Chakraborty [a *] and Dhrubajyoti Chakravarty [b]

[a] Department of Statistics, Dibrugarh University, Dibrugarh, Assam, India.

[b] Department of Statistics, G. C. College, Silchar, Assam, India.

*email: subrata_arya@yahoo.co.in

## Abstract

A new discrete distribution has been proposed as a discrete analogue of the two sided power distribution [*Van Drop, J. R. and Kotz, S. (2002a). A novel extension of the triangular distribution and its parameter estimation, Journal of the Royal Statistical Society, Series D (The Statistician)*, 51, 1: 63-79]. This probability mass function and hazard rate function of this distribution can assume a variety of shapes including bath tub, rectangular, trapezoidal, triangular, J, inverse J, U inverse U, strictly decreasing and strictly increasing shapes. Its moment and reliability properties along with parameter estimation have been investigated.

**Keywords: Two sided power distribution; Hazard rate function; Bathtub shape, Harmonic number**

## 1. INTRODUCTION

The two-sided power (TSP) distribution was first introduced by Van Drop and Kotz (2002a), as an alternative of the beta distribution and an extension of three-parameter triangular distribution which allow J-shaped and U-shaped forms (which are not for the triangular distribution). As the TSP distribution extends the triangular distribution, it inherits its intuitive appeal and parameters interpretation. Van Drop and Kotz (2002a) have dealt with estimation of the parameters of the TSP distribution. Some properties of the two-parameter TSP distribution with support |0, 1| have been discussed in Van Drop and Kotz (2002b).

A random variable $X$ is said to follow a TSP distribution, with parameters $(a, m, b, n)$, $a \le m \le b$, $n > 0$ if its probability density function is given by

$$f(x \mid a, m, b, n) = \begin{cases} \dfrac{n}{b-a}\left(\dfrac{x-a}{m-a}\right)^{n-1} & a \le x \le m \\ \dfrac{n}{b-a}\left(\dfrac{b-x}{b-m}\right)^{n-1} & m \le x \le b \end{cases} \tag{1}$$

It is denoted by $\text{TSP}(a, m, b, n)$. For $n > 1$, the mode of the density function is at '$m$' and the value of the pdf at the mode is always '$n/(b-a)$'. For $0 \le n < 1$ and $a < m < b$ the mode of the density function is at $a$ or $b$ and $f(\bullet \mid a,m,b,n) \to \infty$ at its modes. For $n=1$, $f(\bullet \mid a,m,b,n)$ simplifies to $\text{Uniform}[a,b]$. For $n=2$, $f(\bullet \mid a,m,b,n)$ reduces to a triangular distribution $(a,m,b)$. Finally, for $a=0$ and $m=b=1$, $f(\bullet \mid a,m,b,n)$ corresponds to a power function distribution and for $a=m=0$ and $b=1$ to its reflection.

The cdf of $\text{TSP}(a,m,b,n)$ distribution follows from expression (1) as

$$F(x \mid a,m,b,n) = \begin{cases} \dfrac{m-a}{b-a}\left(\dfrac{x-a}{m-a}\right)^n & a \le x \le m \\ 1-\dfrac{b-m}{b-a}\left(\dfrac{b-x}{b-m}\right)^n & m \le x \le b \end{cases} \quad (2)$$

The survival function of $\text{TSP}(a,m,b,n)$ distribution follows from expression (2) as

$$S(x) = 1 - F(x \mid a,m,b,n) = \begin{cases} 1-\dfrac{m-a}{b-a}\left(\dfrac{x-a}{m-a}\right)^n & a \le x \le m \\ \dfrac{b-m}{b-a}\left(\dfrac{b-x}{b-m}\right)^n & m \le x \le b \end{cases} \quad (3)$$

The mean and variance of the distribution are given respectively by

$$E(X) = \frac{a+(n-1)m+b}{n+1} \quad (4)$$

and $$Var(X) = (b-a)^2 \, \frac{n-2(n-1)(m-a)/(b-a) \times (b-m)/(b-a)}{(n+2)(n+1)^2} \quad (5)$$

**1.1 Interpretations of the parameters**

The parameters $a$ and $b$ are the end point of the support of the distribution, $n$ is the shape parameter and $m$ is the threshold parameter for a change in the form of the pdf. The parameters $a$ and $b$ may be related to pessimistic and optimistic estimates of the associated $\text{TSP}(a,m,b,n)$ random variable.

Given a continuous random variable $X$ with survival function $S(x)$, the discrete analogue (Roy, 2004) is defined as the new random variable $Y = \lfloor X \rfloor$ = largest integer less or equal to *X*. The probability mass function (pmf) of *Y* is then given by

$$\Pr(Y = y) = \Pr(y \le X < y+1) = \Pr(X \ge y) - \Pr(X \ge y+1) = S_X(y) - S_X(y+1) \quad (6)$$

For a detail list of discrete distribution derived from continuous distributions see Chakraborty and Chakravarty (2014).

## 2. DISCRETE TWO-SIDED POWER DISTRIBUTION

The probability mass function (pmf) of a discrete two sided power (DTSP) distribution derived from the pdf of the $\mathrm{TSP}(a,m,b,n)$ distribution using the general discretization approach in equation (6) is given by

$$p_y(a,m,b,n)=\Pr(Y=y)=\begin{cases}\dfrac{(y-a+1)^n-(y-a)^n}{(b-a)(m-a)^{n-1}}, & y=a,a+1,\ldots,m-1\\ \dfrac{(b-y)^n-(b-y-1)^n}{(b-a)(b-m)^{n-1}}, & y=m,m+1,\ldots,b-1\end{cases} \qquad (7)$$

where $a$ , $b$ and $a\le m\le b$ are integers, and $n$ is any positive real number. The random variable $Y$ is said to follow the $\mathrm{DTSP}(a,m,b,n)$. The pmf in the equation (6) is a proper pmf since

$$\sum_{y=a}^{b-1}p_y(a,m,b,n)=\sum_{y=a}^{m-1}\frac{(y-a+1)^n-(y-a)^n}{(b-a)(m-a)^{n-1}}+\sum_{y=m}^{b-1}\frac{(b-y)^n-(b-y-1)^n}{(b-a)(b-m)^{n-1}}$$

$$=\frac{[(1^n-0^n)+(2^n-1^n)+(3^n-2^n)+\cdots+((m-a)^n-(m-a-1)^n)]}{(b-a)(m-a)^{n-1}}$$

$$+\frac{[((b-m)^n-(b-m-1)^n)+(b-m-1)^n-(b-m-2)^n)+\cdots+(2^n-1^n)+(1^n-0^n)]}{(b-a)(b-m)^{n-1}}$$

$$=\frac{(m-a)^n}{(b-a)(m-a)^{n-1}}+\frac{(b-m)^n}{(b-a)(b-m)^{n-1}}$$

$$=\frac{(m-a)^{n-1}(b-m)^{n-1}(m-a+b-m)}{(b-a)(m-a)^{n-1}(b-m)^{n-1}}=1.$$

**Particular cases**:

- When $n=1$, the pmf in equation (7) reduces to discrete *uniform distribution* with pmf

$$p_y(a,m,b,1)=1/(b-a),\ y=a,a+1,\ldots,b-1$$

- For $a<0, b=-a$

$$p_y(a,m,-a,n)=\Pr(Y=y)=\begin{cases}\dfrac{(y-a+1)^n-(y-a)^n}{(-2a)(m-a)^{n-1}}, & y=a,a+1,\ldots,m-1\\ \dfrac{(y+a)^n-(y+a+1)^n}{(2a)(a+m)^{n-1}}, & y=m,m+1,\ldots,-a-1\end{cases}$$

$$p_y(a,-m,-a,n)=\Pr(Y=y)=\begin{cases}\dfrac{(y-a+1)^n-(y-a)^n}{(-1)^n(2a)(a+m)^{n-1}}, & y=a,a+1,\ldots,l-1\\ \dfrac{(y+a)^n-(y+a+1)^n}{(-1)^{n-1}(2a)(m-a)^{n-1}}, & y=l,l+1,\ldots,-a-1\end{cases}$$

Clearly then $p_y(a, m-a, n) = p_{-y-1}(a, -m, -a, n)$ [see corresponding cells of the two rows of figure 1c].

## 2.1 DISTRIBUTIONAL PROPERTIES

### 2.1.1 Pmf plots

The plots of pmf for various values of the parameters have been presented below. It has bee observed that $\mathrm{DTSP}(a, m, b, n)$ can assume a variety of shapes including bath tub, rectangular, trapezoidal, triangular, J, inverse J, strictly decreasing and strictly increasing.

From figure 1a that is when $0 \le n < 1$ and $a < m < b$ the mode of the pmf is at either $a$ or $b$-1 and $p_y \to 0.5$ at its mode while from the figure1b that is for $n > 1$ it has been seen that, the mode of the pmf is at $m$ or $m$-1 the value of the pmf at the mode is always '$n/(b-a+1)$'. From the first plot in figure 1b it can be checked that for $n = 1$, $p_y$ assumes a uniform shape in $[a, b-1]$, for $n = 2$, $p_y$ reduces to a triangular distribution $(a, m, b-1)$ and finally, the plots corresponding cells in the two rows of figure 1c are reflections of each other.

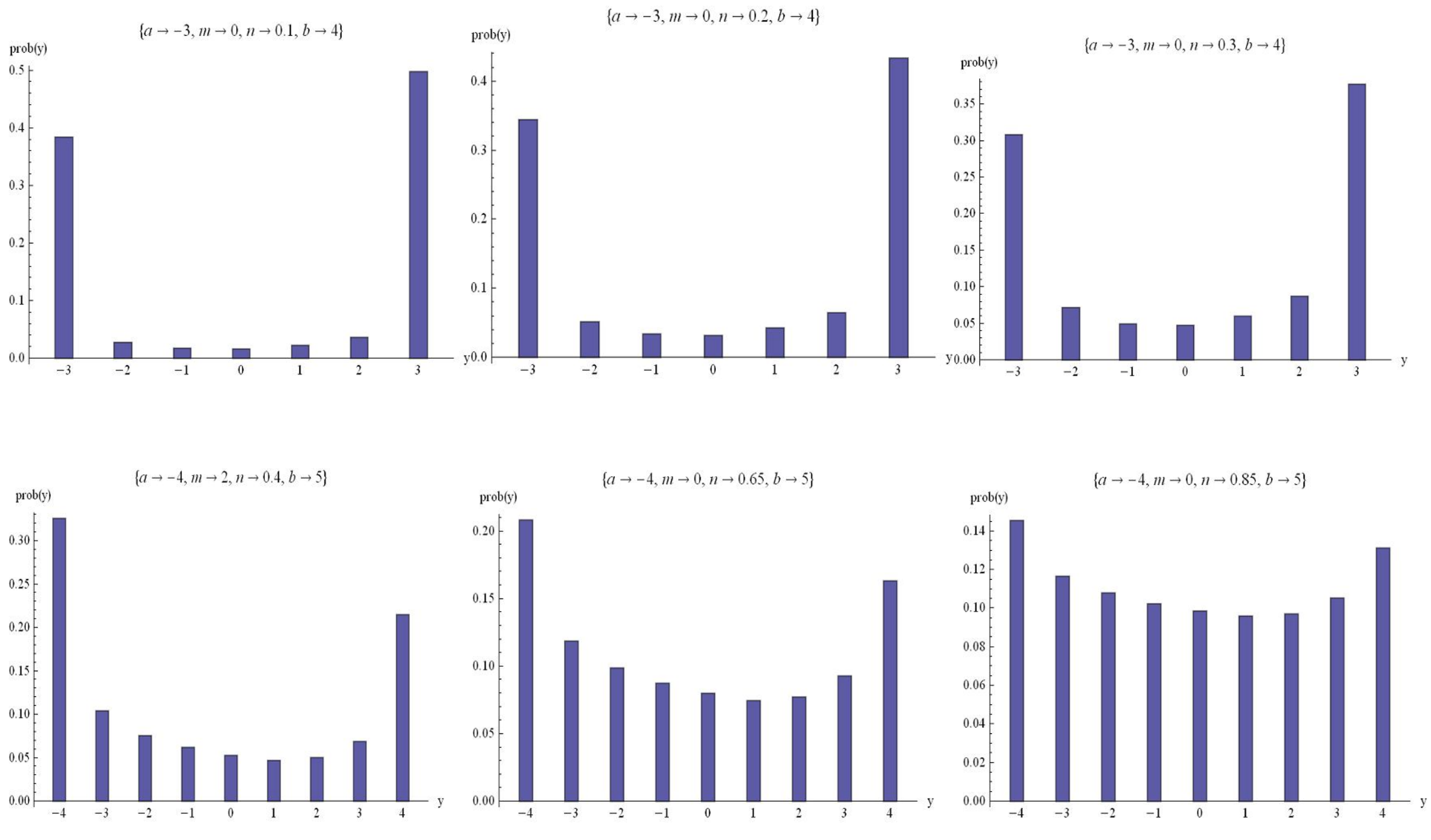

**Figure 1a. pmf of** $\mathrm{DTSP}(a, m, b, n), 0 < n < 1$

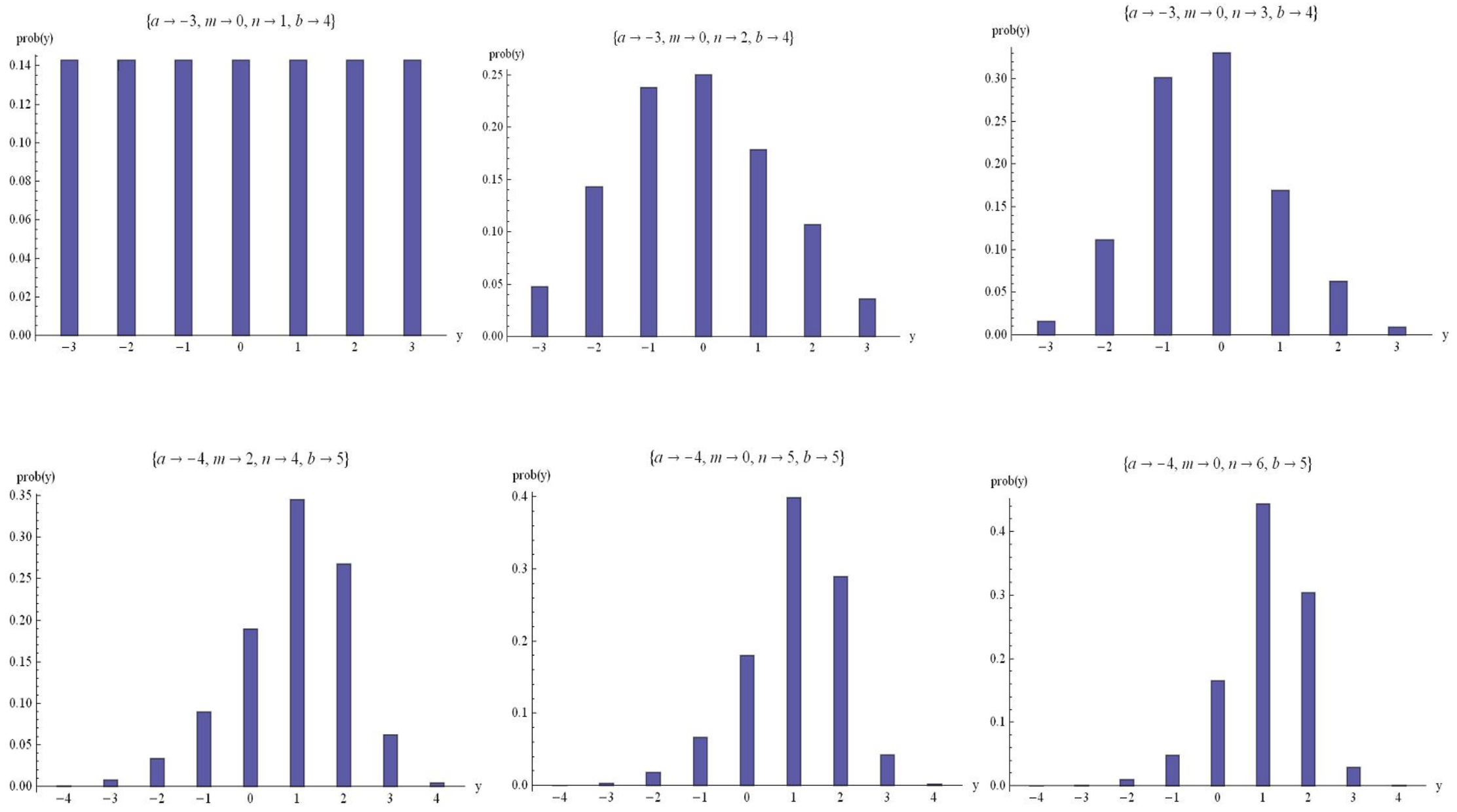

**Figure 1b. pmf of** $\mathrm{DTSP}(a, m, b, n),\ n > 1$

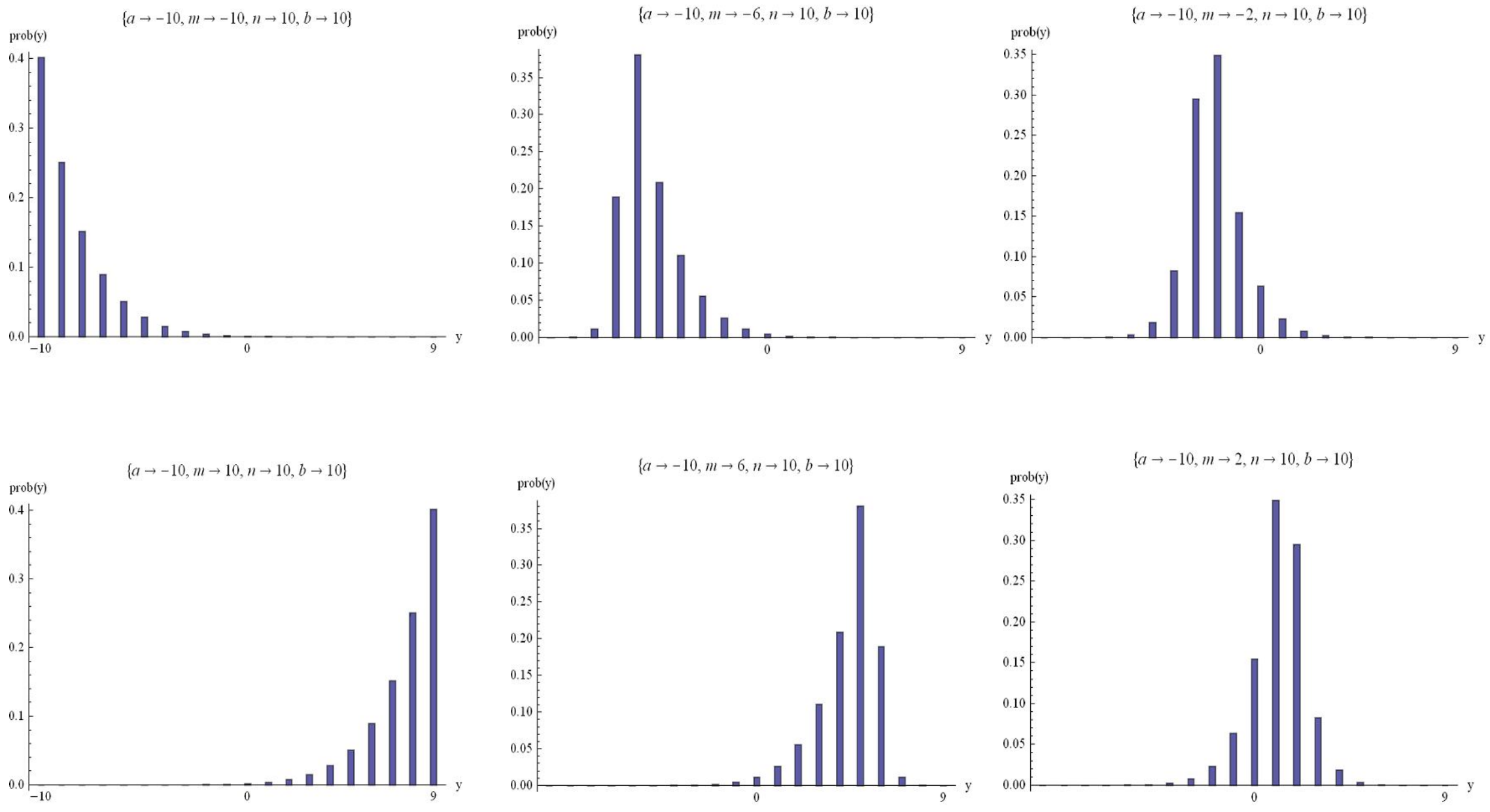

**Figure 1c. pmf of** $\mathrm{DTSP}(a, m, b, n),\ a < 0, b = -a$

### 2.1.2 Recurrence relation for probabilities

The recurrence relation for probabilities of $\mathrm{DTSP}(a, m, b, n)$ is given by

$$p_{y+1}(a,m,b,n)=\begin{cases}\dfrac{(y-a+2)^n-(y-a+1)^n}{(y-a+1)^n-(y-a)^n}\,p_y\,(a,m,b,n),\; y=a,a+1,\ldots,m-1\\[2ex] \dfrac{(b-y-1)^n-(b-y-2)^n}{(b-y)^n-(b-y-1)^n}\,p_y\,(a,m,b,n),\; y=m,m+1,\ldots,b-1\end{cases}$$

### 2.1.3 Cumulative distribution function

**Theorem 1.** The cumulative distribution function (cdf) of $\mathrm{DTSP}(a,m,b,n)$ is given by

$$F(y)=\Pr(Y\le y)=\begin{cases}\dfrac{(m-a)}{(b-a)}\left(\dfrac{y+1-a}{m-a}\right)^n & ,y=a,a+1,\cdots,m-1\\[2ex] 1-\dfrac{(b-m)}{(b-a)}\left(\dfrac{b-y-1}{b-m}\right)^n & ,y=m,m+1,\cdots,b-1\end{cases}$$

### 2.1.4 Moments

The mean and the second moment of $\mathrm{DTSP}(a,m,b,n)$ is respectively

$$E(Y)=\sum_{y=a}^{m-1}y\frac{(y-a+1)^n-(y-a)^n}{(b-a)(m-a)^{n-1}}+\sum_{y=m}^{b-1}y\frac{(b-y)^n-(b-y-1)^n}{(b-a)(b-m)^{n-1}}$$

$$=\frac{\begin{array}{l}[a(1^n-0^n)+(a+1)(2^n-1^n)+(a+2)(3^n-2^n)+\cdots\\ \quad +(m-2)\{(m-a-1)^n-(m-a-2)^n\}+(m-1)\{(m-a)^n-(m-a-1)^n\}]\end{array}}{(b-a)(m-a)^{n-1}}$$

$$+\frac{\begin{array}{l}[m\{(b-m)^n-(b-m-1)^n\}+(m+1)\{(b-m-1)^n-(b-m-2)^n\}+\cdots\\ \quad +(b-2)(2^n-1^n)+(b-1)(1^n-0^n)\end{array}}{(b-a)(b-m)^{n-1}}$$

$$=\frac{(m-1)(m-a)^n-(1^n+2^n+3^n+\cdots+(m-a-1)^n)}{(b-a)(m-a)^{n-1}}+\frac{m(b-m)^n-(1^n+2^n+3^n+\cdots+(b-m-1)^n)}{(b-a)(b-m)^{n-1}}$$

$$=\frac{(m-a)^{1-n}H_{(m-a-1)}^{(-n)}+a\,(m-1)+(b-m)^{1-n}H_{(b-m-1)}^{(-n)}-m(b-1)}{(a-b)} \qquad (8)$$

and $E(Y^2)=\sum_{y=a}^{m-1}y^2\dfrac{(y-a+1)^n-(y-a)^n}{(b-a)(m-a)^{n-1}}+\sum_{y=m}^{b-1}y^2\dfrac{(b-y)^n-(b-y-1)^n}{(b-a)(b-m)^{n-1}}$

$$
= \frac{\begin{array}{l}[a^2(1^n-0^n)+(a+1)^2(2^n-1^n)+(a+2)^2(3^n-2^n)+\cdots \\ \qquad +(m-2)^2\{(m-a-1)^n-(m-a-2)^n\}+(m-1)^2\{(m-a)^n-(m-a-1)^n\}]\end{array}}{(b-a)(m-a)^{n-1}}
$$

$$
+ \frac{\begin{array}{l}[m^2\{(b-m)^n-(b-m-1)^n\}+(m+1)^2\{(b-m-1)^n-(b-m-2)^n\}+\cdots \\ \qquad +(b-2)^2(2^n-1^n)+(b-1)^2(1^n-0^n)\end{array}}{(b-a)(b-m)^{n-1}}
$$

$$
= \frac{(m-1)^2(m-a)^n+\{a^2-(a+1)^2\}1^n+\{(a+1)^2-(a+2)^2\}2^n+\cdots+\{(m-2)^2-(m-1)^2\}(m-a-1)^n}{(b-a)(m-a)^{n-1}}
$$

$$
+ \frac{m(b-m)^n+\{(b-1)^2-(b-2)^2\}1^n+\{(b-1)^2-(b-2)^2\}2^n+\cdots+\{(m+1)^2-m^2\}(b-m-1)^n}{(b-a)(b-m)^{n-1}}
$$

$$
= \frac{(m-1)^2(m-a)^n+(-2a-1)1^n+(-2a-3)2^n+\cdots+(-2m+3)(m-a-1)^n}{(b-a)(m-a)^{n-1}}
$$

$$
+ \frac{m(b-m)^n+(2b-3)1^n+(2b-5)2^n+\cdots+(2m+1)(b-m-1)^n}{(b-a)(b-m)^{n-1}}
$$

$$
= \frac{(m-1)^2(m-a)^n+\sum_{i=1}^{m-a-1}[-2a-(2i-1)]i^n}{(b-a)(m-a)^{n-1}} + \frac{m(b-m)^n+\sum_{i=1}^{b-m-1}[2b-(2i+1)]i^n}{(b-a)(b-m)^{n-1}}
$$

$$
= \frac{(m-1)^2(m-a)^n-(2a-1)H_{(m-a-1)}^{(-n)}-2H_{(m-a-1)}^{(-n-1)}}{(b-a)(m-a)^{n-1}} + \frac{m(b-m)^n+(2b-1)H_{(b-m-1)}^{(-n)}-2H_{(b-m-1)}^{(-n-1)}}{(b-a)(b-m)^{n-1}} \tag{9}
$$

Hence the variance is

$$
Var(Y) = \frac{(m-1)^2(m-a)^n-(2a-1)H_{(m-a-1)}^{(-n)}-2H_{(m-a-1)}^{(-n-1)}}{(b-a)(m-a)^{n-1}} + \frac{m(b-m)^n+(2b-1)H_{(b-m-1)}^{(-n)}-2H_{(b-m-1)}^{(-n-1)}}{(b-a)(b-m)^{n-1}}
$$
$$
+ \left( \frac{(m-a)^{1-n}H_{(m-a-1)}^{(-n)}+a(m-1)+(b-m)^{1-n}H_{(b-m-1)}^{(-n)}-m(b-1)}{(a-b)} \right)^2,
$$

where $H_{(a)}^{(b)} = \sum_{k=1}^{a}(1/k^b)$ is the generalized harmonic number of order $b$, whose limit as $a \to \infty$ exists when $b > 1$. The related sum $\sum_{k=1}^{a} k^b$ occurs in the study of Bernoulli numbers, the harmonic number also appears in the study of Stirling's numbers (Abramowitz and Stegun, 1970).

**Theorem 2.** The mean and variance of $\mathrm{DTSP}(a,m,b,n)$ distribution are bounded.

***Proof***. The discretized version $Y$ of $X$ that is $\mathrm{DTSP}(a,m,b,n)$ is defined as $Y = \lfloor X \rfloor$ =largest integer less or equal to $X$ . It can be assumed that $X = Y + U$ , where $U$ is the fractional part of $X$ which is chopped off from $X$ to obtain $Y$. Therefore, $U$ will have another continuous probability distribution in the support $(0,1)$, which is independent of $Y$. Hence $E(Y)=E(X-U)=E(X)-E(U)$. But $0<E(U)<1$, therefore on using equation (4)

$$\frac{a+(n-1)m+b}{n+1}-1<E(U)<\frac{a+(n-1)m+b}{n+1}$$

Similar argument gives

$V(Y)=V(X-U)=V(X)+V(U) \quad \Rightarrow V(Y)=V(X)+V(U)$ [Assuming independence of $X$ and $U$].

But for any continuous random variable $U$ in $(0,1)$, $0<Var(U)\le 1/4$, now on using equation (5)

$$(b-a)^2\frac{n-2(n-1)(m-a)/(b-a)\times(b-m)/(b-a)}{(n+2)(n+1)^2}< V(Y)$$

$$<(b-a)^2\frac{n-2(n-1)(m-a)/(b-a)\times(b-m)/(b-a)}{(n+2)(n+1)^2}+0.25$$

Hence the mean and variance is bounded above (see also figures 2 and 4).

Plots of the mean, variance and ID for different combination of parameter values for $\mathrm{TSP}(a,m,b,n)$ and $\mathrm{DTSP}(a,m,b,n)$ have been presented below.

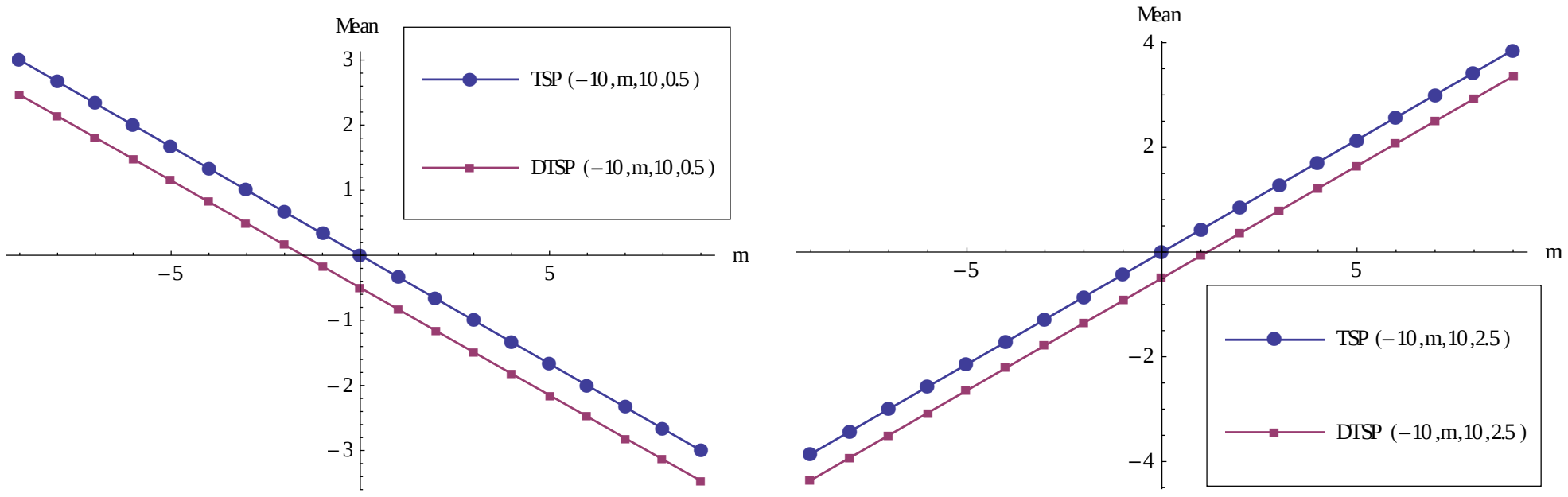

**Figure 2. Mean of $\mathrm{TSP}(a,m,b,n)$ and $\mathrm{DTSP}(a,m,b,n)$ for varying *m***

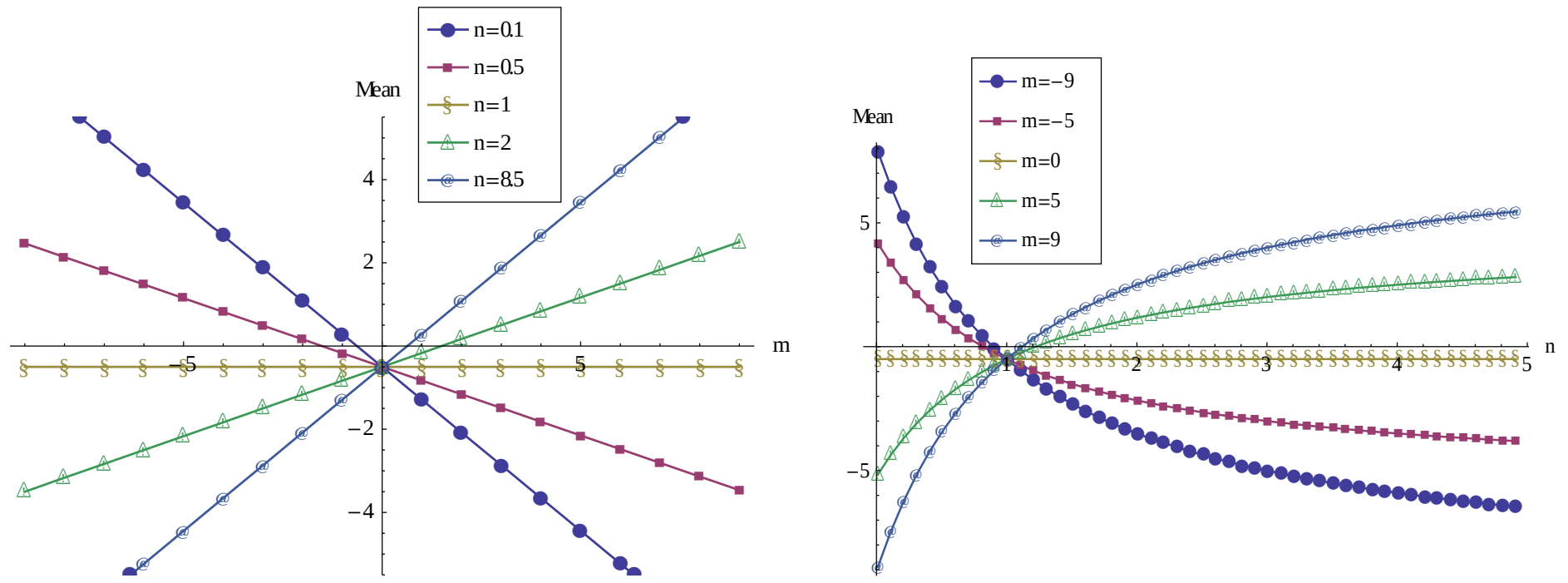

**Figure 3. Mean of $\mathrm{DTSP}(-10,m,10,n)$ for varying *m* and *n***

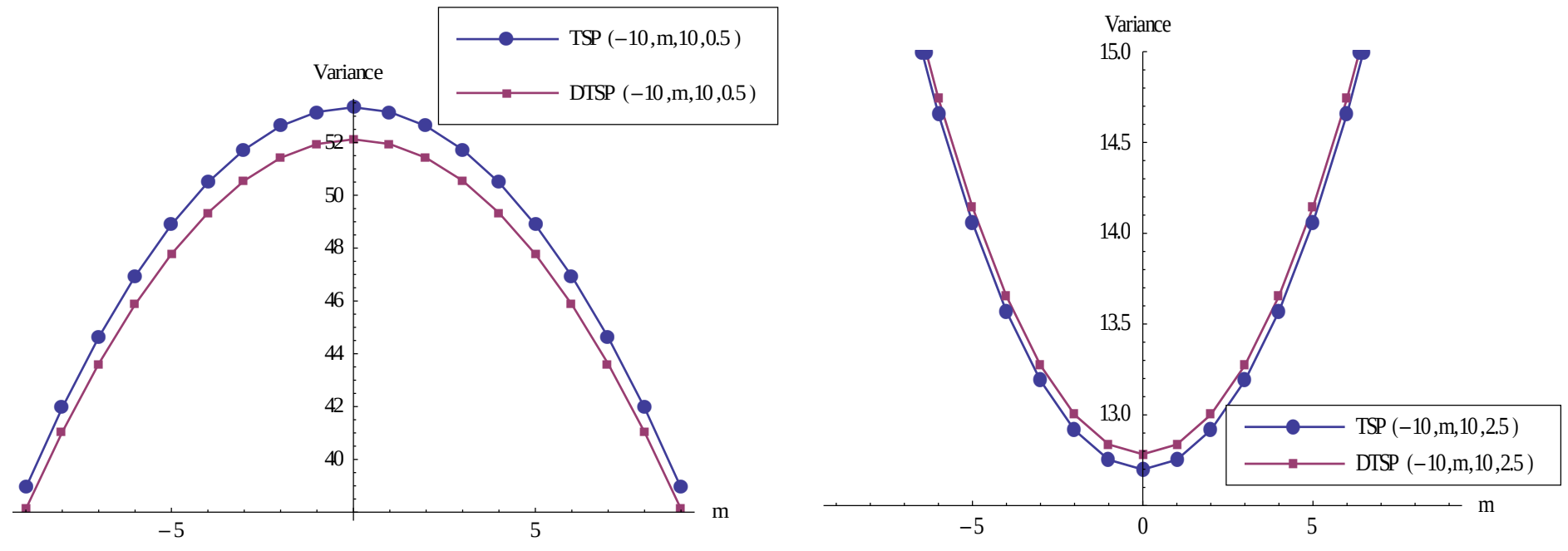

**Figure 4. Variance of $\mathrm{TSP}(a,m,b,n)$ and $\mathrm{DTSP}(a,m,b,n)$ for varying *m***

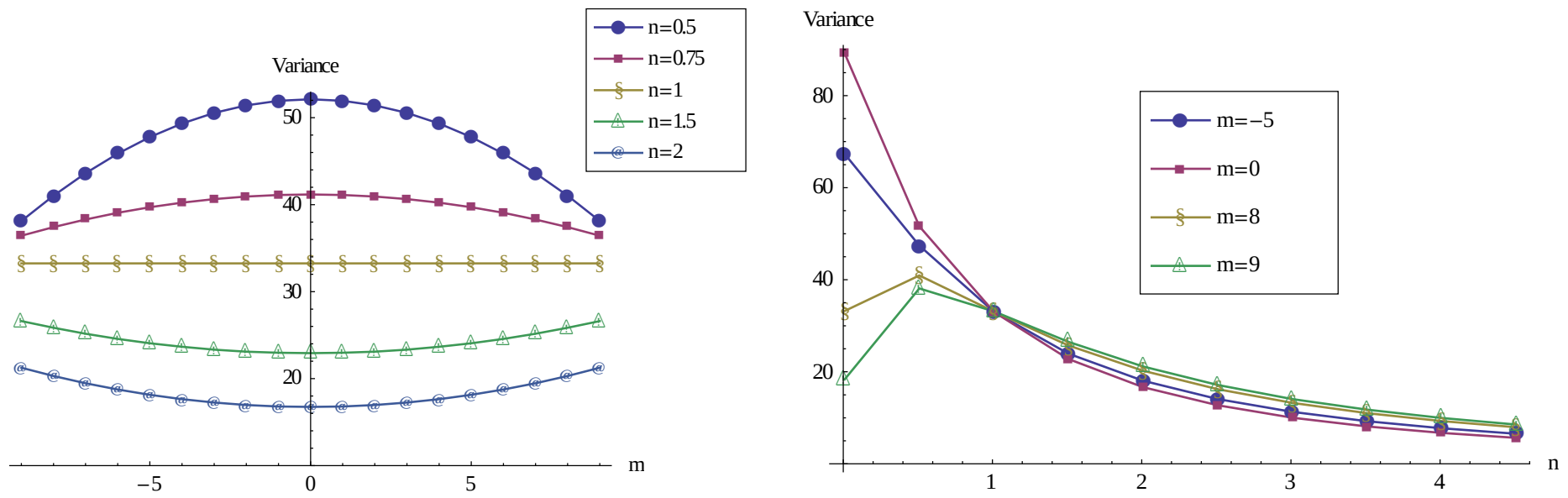

**Figure 5. Variance of $\mathrm{DTSP}(-10,m,10,n)$ for varying *m* and *n***

### 2.1.5 Mode

For $n>1$, the mode of the pmf is either at $m$ or at $m-1$. For $0\le n<1$ and $a<m<b$, the mode of the pmf is always either at '$a$' or at '$b-1$' with respective value of the pmf equal to $1/((b-a)(b-m)^{n-1})$ and $1/((b-a)(m-a)^{n-1})$ as its mode. This can also be seen from figure 1a above.

## 2.2 RELIABILITY PROPERTIES

### 2.2.1 Survival function

**Theorem 3.** The survival function of $\mathrm{DTSP}(a,m,b,n)$ is given by

$$S(y)=\begin{cases}1-\dfrac{(m-a)}{(b-a)}\left(\dfrac{y-a}{m-a}\right)^{n} & ,y=a,a+1,\cdots,m-1\\[2ex] \dfrac{(b-m)}{(b-a)}\left(\dfrac{b-y}{b-m}\right)^{n} & ,y=m,m+1,\cdots,b-1\end{cases}\tag{10}$$

***Proof.*** The survival function is defined by $S(y)=P(Y\ge y)$.

**Case I.** $y\le m-1$

$$S(y)=[p(y)+p(y+1)+\cdots+p(m-1)]+[p(m)+p(m+1)+\cdots+p(b-2)+(b-1)]$$

$$=\left[\frac{(y-a+1)^n-(y-a)^n}{(b-a)(m-a)^{n-1}}+\frac{(y-a+2)^n-(y-a+1)^n}{(b-a)(m-a)^{n-1}}+\cdots+\frac{(m-a)^n-(m-a-1)^n}{(b-a)(m-a)^{n-1}}\right]$$

$$+\left[\frac{(b-m)^n-(b-m-1)^n}{(b-a)(b-m)^{n-1}}+\frac{(b-m-1)^n-(b-m-2)^n}{(b-a)(b-m)^{n-1}}+\cdots+\frac{1^n-0^n}{(b-a)(b-m)^{n-1}}\right]$$

$$=\frac{(m-a)^n-(y-a)^n}{(b-a)(m-a)^{n-1}}+\frac{b-m}{b-a}=1-\frac{(m-a)}{(b-a)}\left(\frac{y-a}{m-a}\right)^n$$

**Case II.** $y \geq m$

$$S(y)=[p(y)+p(y+1)+\cdots+p(b-1)]$$

$$=\left[\frac{(b-y)^n-(b-y-1)^n}{(b-a)(b-m)^{n-1}}+\frac{(b-y-1)^n-(b-y-2)^n}{(b-a)(b-m)^{n-1}}+\cdots+\frac{1^n-0^n}{(b-a)(b-m)^{n-1}}\right]$$

$$=\frac{(b-y)^n}{(b-a)(b-m)^{n-1}}=\frac{(b-m)}{(b-a)}\left(\frac{b-y}{b-m}\right)^n$$

Combining these two cases gives the desires result.

The plots of survival function for various pair of parametric values have been presented in figure 7 below.

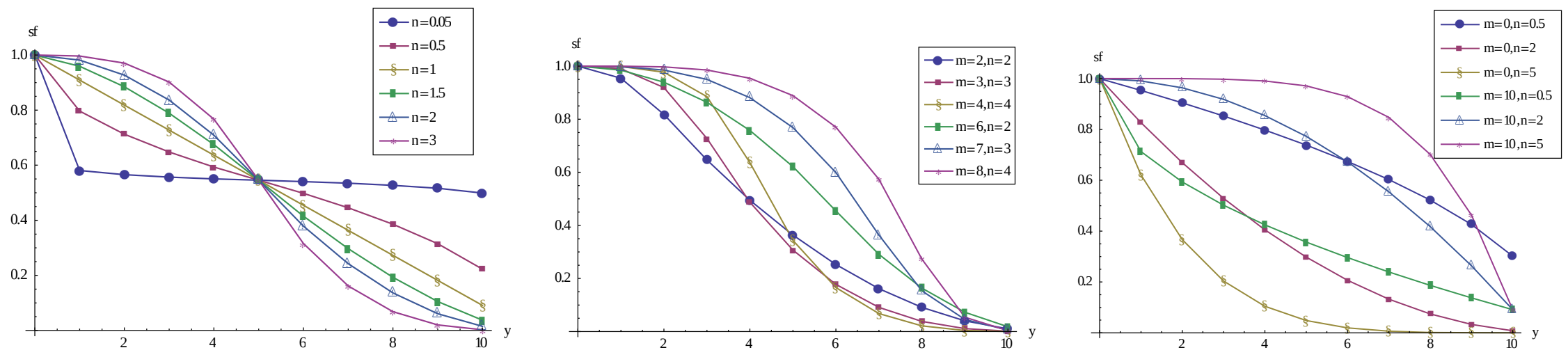

**Figure 7. Survival function of** $\mathrm{DTSP}(0,5,11,n)$**,** $\mathrm{DTSP}(0,m,11,n)$ **and** $\mathrm{DTSP}(0,m,11,n)$

***Remark* 1**. It is to be noted that the survival function of $\mathrm{DTSP}(a,m,b,n)$ is same as that of continuous two sided power distribution given in equation (3).

***Remark* 2**. From equation (10), $S(m)=(b-m)/(b-a)$. Therefore $m$ is the $(1-(b-m)/(b-a))^{th}$ or simply $((m-a)/(b-a))^{th}$ quantile of the distribution irrespective of the choice of the parameter $n$. In fact $m$ will be the median if it is equal to $\lfloor (a+b)/2 \rfloor$, since

$$S(m)=1/2 \Rightarrow (b-m)/(b-a)=1/2 \Rightarrow m=(a+b)/2$$

Thus the parameter '$m$' identifies $100\,((m-a)/(b-a))^{th}$ percentile or simply $((m-a)/(b-a))^{th}$ quantile of $\mathrm{DTSP}(a,m,b,n)$ distribution.

### 2.2.2 Hazard (Failure) rate function

From equation (10), it is a simple exercise to obtain the failure or hazard rate function of a discrete random variable following the $\mathrm{DTSP}(a,m,b,n)$ distribution with pmf (7), given by

$$r(y)=\frac{\Pr(Y=y)}{S(y)}=\begin{cases}\dfrac{(y-a+1)^n-(y-a)^n}{(b-a)(m-a)^{n-1}-(y-a)^n}, & y=a,a+1,\ldots,m-1\\[2ex] \dfrac{(b-y)^n-(b-y-1)^n}{(b-y)^n}, & y=m,m+1,\ldots,b-1\end{cases}$$

The hazard rate function for different pairs of values of parameters has been graphed in figure 8

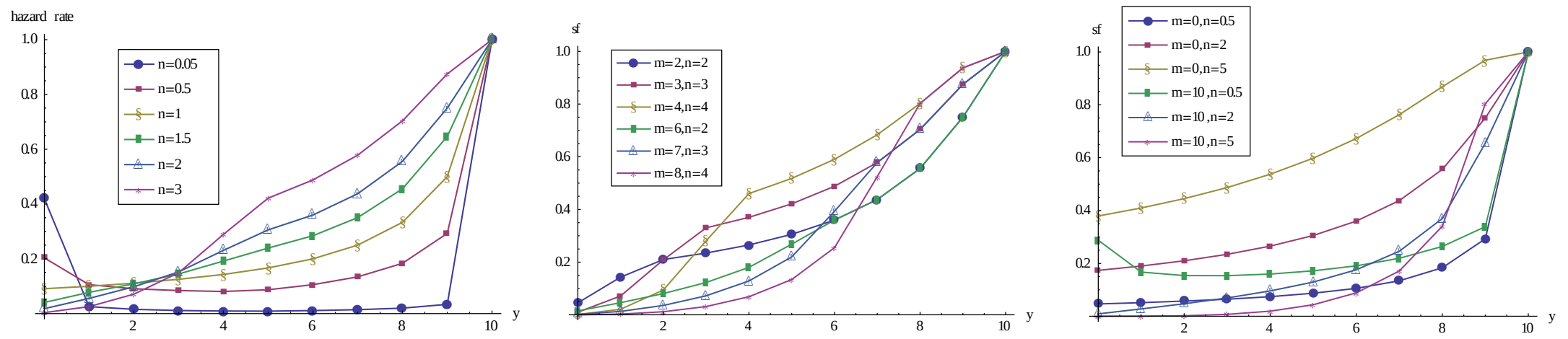

**Figure 8. Hazard rate function plot of** $\mathrm{DTSP}(0,5,11,n)$**,** $\mathrm{DTSP}(0,m,11,n)$**,** $\mathrm{DTSP}(0,m,11,n)$

From the above graphs, it can be seen that $\mathrm{DTSP}(a,m,b,n)$ has most of hazard rate shapes including the Bathtub shaped one, except possibly the strictly decreasing one.

## 2.3 Estimation of parameters

For $\mathrm{DTSP}(a,m,b,n)$, the value of $a$ and $b$ can be obtained from the data given, as these are two end points of the distribution. Also $m$ is the threshold parameter which indicated the change point of the pmf and can be realistically estimated from the data. Therefore the only need is to estimate the parameter $n$. This parameter can be estimated by different methods as discussed below

### 2.3.1 Maximum likelihood estimation (MLE)

For a random sample $y_1, y_2, \cdots, y_k$ of size $k$ the log likelihood function of $\mathrm{DTSP}(a,m,b,n)$ is given by

$$\log L=\sum_{i=a}^{m-1}\log\left(\frac{(y_i-a+1)^n-(y_i-a)^n}{(b-a)(m-a)^{n-1}}\right)+\sum_{i=m}^{b-1}\log\left(\frac{(b-y_i)^n-(b-y_i-1)^n}{(b-a)(b-m)^{n-1}}\right)$$

Then the likelihood equation is given by

$$\log(b-m)+\log(m-a)-\sum_{i=a}^{m-1}\frac{(y_i-a+1)^n\log(y_i-a+1)-(y_i-a)^n\log(y_i-a)}{(y_i-a+1)^n-(y_i-a)^n}$$
$$-\sum_{i=m}^{b-1}\frac{(b-y_i)^n\log(b-y_i)-(b-y_i-1)^n\log(b-y_i-1)}{(b-y_i)^n-(b-y_i-1)^n}=0$$

The MLE of $n$ can be obtained by numerically solving the likelihood equation or by maximizing the log-likelihood function by global optimization method.

### 2.3.2 Method of moment estimation (MME)

The unknown parameter $n$ can be estimated by solving the equation $E(Y)-m_1=0$, where $m_1$ is the sample mean. Alternatively, $n$ can be estimated by minimizing $[E(Y)-m_1]^2$ with respect to $n$ (Khan et al. 1989).

### 2.4 Simulation study for MLE and MMEs

A through simulation study has been carried out to ascertain the effectiveness of the parameter estimates of the $\text{DTSP}(a,m,b,n)$ by generating 1000 samples of different sizes for different choices of the parameter $n$. Here the random integer ($Y$) from $\text{DTSP}(a,m,b,n)$ has been sampled by first using the inverse transformation method to generate from continuous distribution and then taking the floor of the values of the resulting continuous variate. The steps for generating a random integer are follows:

Step I: Consider a continuous random variable $X$ with distribution function $F(x)$.

Step II: Generate a random number $U$ from uniform distribution $U(0,1)$.

Step III: Compute $X = F^{-1}(U)$ .

Step IV: $Y = \lfloor X \rfloor$.

The result of simulation analysis has been presented in tables 1 and 2. All the entries of the table are the means of estimates of 1000 samples. Accuracy and precision of the estimation methods have been checked and established using the estimates of the following criteria:

I. Estimate of $E(\hat{\theta})$ by $(1/k)\sum_{i=1}^{k}\hat{\theta}_i$

II. Estimate of bias of $\hat{\theta}$: $\text{Bias} = (1/k)\sum_{i=1}^{k}(\hat{\theta}_i-\theta)$

III. Estimate of mean square error (MSE): $\mathrm{MSE}(\hat{\theta}) = (1/k)\sum_{i=1}^{k}(\hat{\theta}_i - \theta)^2$

IV. Estimate of Variance: $\mathrm{Var}(\hat{\theta}) = (1/k)\sum_{i=1}^{k}\left(\hat{\theta}_i^2 - E(\hat{\theta})\right)^2$

where $\hat{\theta}_i$ is the estimate of the unknown true value $\theta$ obtained from the $i^{\text{th}}$ sample $(i = 1,2,\cdots,k)$. Here $k$ = 1000.

In addition, to check the asymptotic normality of the estimators the 95% confidence intervals for the estimators using the formula $\mathrm{CI}_{95\%} = \hat{\theta} \pm 1.96\,\mathrm{SE}(\hat{\theta})$, where $\mathrm{SE}(\hat{\theta}) = \sqrt{\mathrm{Var}(\hat{\theta})}$ have been computed and the percentage of estimates $\hat{\theta}_i$ belonging to the interval obtained.

**Table 1. Results of simulation from $\mathrm{DTSP}(a, m, b, n)$ distribution for $a = -10, m = 0, b = 10, n = 0.5$**

| Parameter value | Sample size→ | 25 | | 50 | | 100 | |
|---|---|---|---|---|---|---|---|
| | estimates↓ | MLE | MME | MLE | MME | MLE | MME |
| $n = 0.5$ | $E(\hat{n})$ | 0.6087 | 0.4224 | 0.5950 | 0.4095 | 0.5868 | 0.4052 |
| | Bias( $\hat{n}$ ) | -0.1087 | 0.0776 | -0.0950 | 0.0905 | -0.0868 | 0.0948 |
| | $MSE(\hat{n})$ | 0.0360 | 0.0190 | 0.0197 | 0.0133 | 0.0131 | 0.0114 |
| | $Var(\hat{n})$ | 0.0242 | 0.0130 | 0.0107 | 0.0051 | 0.0056 | 0.0024 |
| | % of $n$ in CI | 95.30 | 94.40 | 95.40 | 95.60 | 94.60 | 94.90 |

**Table 2. Results of simulation from $\mathrm{DTSP}(a, m, b, n)$ distribution for $a = -10, m = 0, b = 10, n = 3.5$**

| Parameter value | Sample size→ | 25 | | 50 | | 100 | |
|---|---|---|---|---|---|---|---|
| | estimates↓ | MLE | MME | MLE | MME | MLE | MME |
| $n = 3.5$ | $E(\hat{n})$ | 1.7344 | 3.5534 | 1.9604 | 3.4511 | 2.1736 | 3.3887 |
| | Bias( $\hat{n}$ ) | 1.7656 | -0.0537 | 1.5396 | 0.0489 | 1.3264 | 0.1113 |
| | $MSE(\hat{n})$ | 3.2829 | 0.6981 | 2.4924 | 0.3088 | 1.8553 | 0.1614 |
| | $Var(\hat{n})$ | 0.1656 | 0.6959 | 0.1222 | 0.3067 | 0.0960 | 0.1492 |
| | % of $n$ in CI | 95.70 | 95.70 | 95.50 | 95.80 | 94.90 | 96.00 |

***Remark* 3**. From the tables 1 and 2 it can be seen that larger the sample size, smaller is the MSE. But MME seems to be the better option among the two with respect to the different criteria.

**3. Conclusion**

In this paper a discrete two sided power distribution has been derived. Various distributional and reliability properties have been derived. Mean and variance has been obtained in compact form. Parameter estimation through maximum likelihood and moment method has been considered. It is envisaged that this discrete distribution which can model be suitable for modeling U shaped, J shaped and inverse J shaped, bathtub shaped mass function and variety of failure time count data will be a valuable addition to the repository of discrete distribution. Further works with this distribution will be followed up later.